\documentclass[a4paper]{article}
\usepackage{RR}
\usepackage{hyperref}
\usepackage[frenchb]{babel} 

\RRdate{Septembre 2006}

\RRauthor{
Abdou W. BELLO\thanks[sfn]{Universit\'e d'Abomey-Calavi, 01 BP 526 Cotonou, B\'enin}
  \and
Aur\'elien GOUDJO\thanksref{sfn}
\and
C\^ome GOUDJO\thanksref{sfn}
\and
Herv\'e GUILLARD\thanks[sf]{INRIA, 2004 Route des Lucioles, BP 93, F-06902 Sophia-Antipolis cedex, France}
\and
Jean-Antoine DESIDERI\thanksref{sf}%

}
\authorhead{A.W.BELLO, A.GOUDJO, C.GOUDJO, H.GUILLARD \& J.-A.DESIDERI}
\RRtitle{Un sch\'ema de type Volumes-Finis-Roe (VFRoe) pour les \'equations de Saint-Venant 1D : \\ {\LARGE Simulation num\'erique des bancs-couvrants-d\'ecouvrants}}
\RRetitle{A VFRoe scheme for 1D shallow water flows : \\wetting and drying simulation}
\titlehead{Equations de Saint-Venant : traitement des bancs couvrants-d\'ecouvrants}
\RRresume{Nous pr\'esentons dans ce rapport la r\'esolution, par une m\'ethode volumes finis, du syst\`eme des \'equations de Saint-Venant avec termes sources topographiques sur des domaines 1D. Avec une id\'ee originale de Leroux \cite{Leroux}, le syst\`eme des \'equations est compl\'et\'e par une \'equation triviale sur la bathym\'etrie. Par un changement de variable, on \'elabore une formulation c\'el\'erit\'e-vitesse des \'equations que l'on lin\'earise. Nous construisons ensuite un solveur de Riemann approch\'e qui pr\'eserve la positivit\'e de la c\'el\'erit\'e et qui assure la prise en compte des {\it bancs-couvrants-d\'ecouvrants}. Enfin, des applications num\'eriques sur des cas tests sont pr\'esent\'ees.
}

\RRabstract{A finite-volume method for the one-dimensional shallow-water equations 
including topographic source terms is presented. Exploiting an original 
idea by Leroux \cite{Leroux}, the system of partial-differential equations 
is completed by a trivial equation for the bathymetry. By applying a 
change of variable, the system is given a celerity-speed formulation,
and linearized. As a result, an approximate Riemann
solver preserving the positivity of the celerity can be constructed,
permitting {\it wetting and drying flow} simulations to be performed.
Finally, the simulation of numerical test cases is presented.}
\RRmotcle{Equations de Saint-Venant, volumes finis, solveur de Riemann, sch\'ema positif, bancs couvrants-d\'ecouvrants}
\RRkeyword{Shallow water equations, finite volumes, Riemann solver, positivity preserving scheme, wetting and drying flows}
\RRprojet{Opale}  
\RRtheme{\THNum} 
\URSophia 
\usepackage{amsfonts,hhline,ulem,delarray,amssymb,amsmath}

\usepackage{amssymb}
\def\N{\mathbb{N} }
\def\R{\mathbb{R} }

\newtheorem{prop}{Proposition}[section]

\newtheorem{rmq}{Remarque}[section]

\begin{document}
\makeRR   
\tableofcontents

\newpage
\section{Introduction}
\begin{figure}[ht]
\begin{picture}(0,220)
\includegraphics[scale=.43]{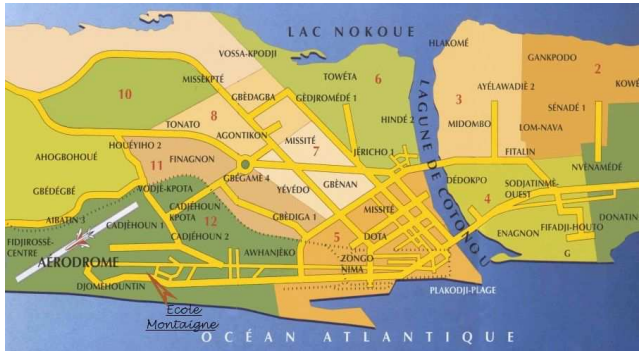}
\end{picture}
\caption{\textbf{Chenal de Cotonou}}
\label{fig1}
\end{figure}

La ville de Cotonou est dans une cuvette ceintur\'ee de plans d'eau et travers\'ee du nord au sud par un chenal (la Lagune de Cotonou) de 4000 m de long sur 350 m de large, reliant le lac Nokou\'e \`a l'oc\'ean Atlantique (Fig.\ref{fig1}).

Les quelques ouvrages d'assainissement dont dispose la ville sont g\'en\'eralement en bordure des routes et, pour la plupart, connect\'es au chenal. Malheureusement, ce syst\`eme de drainage au lieu d'aider \`a l'\'evacuation des eaux pluviales, sert plut\^ot de vecteur \`a l'invasion de la ville par les eaux de crue du chenal.

Les \'ecoulements d'un tel r\'eseau d'adduction d'eau sont des \'ecoulements \`a surface libre, en eaux peu profondes (shallow water) et sont alors r\'egis par un syst\`eme bidimensionnel des \'equations de Saint-Venant. Ces \'equations sont obtenues \`a partir des \'equations de Navier-Stokes pour un fluide incompressible en faisant l'hypoth\`ese de pression hydrostatique, de vitesses uniformes suivant la verticale, d'un fond et d'une surface libre imperm\'eables. On les emploie dans des domaines aussi divers que la protection de l'environnement, le calcul des mar\'ees et des ondes de temp\^ete, la s\'edimentologie, la simulation des ondes de submersion, l'\'etude des crues, etc.

Dans ce rapport, comme premi\`ere \'etape dans la r\'ealisation d'une simulation num\'erique r\'ealiste, nous proposons une \'etude du syst\`eme unidimensionnel des \'equations de Saint-Venant avec topographie. Suivant l'id\'ee d\'evelopp\'ee dans \cite{gallouet, gallouet2, gallouet3}, nous utilisons une m\'ethode de type Volumes-Finis-Roe (VFRoe) se basant sur une formulation c\'el\'erit\'e-vitesse des \'equations. Une telle formulation peut conduire \`a une c\'el\'erit\'e n\'egative de l'\'etat interm\'ediaire du probl\`eme de Riemann y d\'ecoulant. Notre contribution est la construction d'un solveur de Riemann approch\'e  avec un choix des vitesses d'ondes garantissant la pr\'eservation de la positivit\'e de la c\'el\'erit\'e de l'\'etat interm\'ediaire. 

Ensuite, gr\^ace aux r\'esultats de \cite{david, leveque} nous adoptons le solveur de Riemann approch\'e afin de prendre en compte le cas des bancs-couvrants-d\'ecouvrants (zones noy\'ees/non-noy\'ees). Le traitement de ce type de cas est essentiel en vue de la simulation num\'erique des probl\`emes d'intrusion d'eau dans la ville de Cotonou.

Enfin, pour tester la robustesse de notre sch\'ema, dans un premier temps nous reprenons  des cas-tests effectu\'es dans \cite{Leroux} puis dans un second temps deux derniers cas-tests sont pr\'esent\'es pour tester le d\'ebordement de l'eau d'un canal (bancs-couvrants-d\'ecouvrants).

\section{Le mod\`ele math\'ematique}

En absence des forces de frottement, le syst\`eme 1D des \'equations de Saint-Venant est donn\'e par :
\begin{eqnarray}\label{e1}\begin{array}\{{lcl}.
\displaystyle\frac{\partial h}{\partial t}+\frac{\partial (hu)}{\partial x}&=&0\\&&\\\displaystyle\frac{\partial (hu)}{\partial t}+ \frac{\partial }{\partial x}\left(hu^2 + \frac{1}{2}gh^2 \right)&=&\displaystyle-gh\frac{\partial a}{\partial x}\end{array},(x,t)\in\R\times\R_+\end{eqnarray}
$u(x,t)$ est la vitesse de l'eau, $h(x,t)$ la hauteur d'eau, $a(x)$ la hauteur de la topographie du sol, $h+a$ l'\'elevation de la surface libre de l'eau ($a$ et $h+a$ sont prises par rapport \`a un plan de r\'ef\'erence), $g$ l'acc\'el\'eration gravitationnelle.

La bathym\'etrie $a$ \'etant ind\'ependante du temps, on compl\`ete le syst\`eme (\ref{e1}) par l'\'equation $\dfrac{\partial a}{\partial t}=0$ \cite{Leroux} et on obtient :

\begin{eqnarray}\label{e2}\begin{array}\{{lcl}.
\displaystyle\frac{\partial h}{\partial t}+\frac{\partial q}{\partial x}&=&0\\&&\\\displaystyle\frac{\partial q}{\partial t}+ \frac{\partial }{\partial x}\left(\dfrac{q^2}{h} + \frac{1}{2}gh^2 \right)+gh\dfrac{\partial a}{\partial x}&=&0\\&&\\\dfrac{\partial a}{\partial t}&=&0\end{array},(x,t)\in\R\times\R_+,\end{eqnarray} o\`u $q=hu$ repr\'esente le d\'ebit d'eau.

\section{Int\'egration num\'erique}

Nous pr\'esentons dans cette section l'int\'egration des \'equations par utilisation des techniques de volumes finis et de solveurs de Riemann.

\subsection{Le maillage}

On consid\`ere un \'ecoulement unidimensionnel sur une longueur $L$.

L'intervale $[0,L]$ est subdivis\'e en $N$ segments de m\^eme amplitude $\Delta x=\dfrac{L}{N}$. On obtient alors une suite de points $(x_j)_{j\in J=\{0,\ldots,N\}}$ d\'efinis par : $x_j=j\Delta x$. On pose 

\begin{eqnarray}\begin{array}\{{rcl}. x_{j+\frac{1}{2}}&=&\dfrac{1}{2}(x_j+x_{j+1}),\;j\in\{0,1,\ldots,N-1\}\\&&\\C_j&=&\left]x_{j-\frac{1}{2}};x_{j+\frac{1}{2}}\right[,\;j\in\{1,\ldots,N-1\}\\&&\\C_0&=&\left]x_{0};x_{\frac{1}{2}}\right[\\&&\\C_N&=&\left]x_{N-\frac{1}{2}};x_N\right[\end{array}\end{eqnarray}

\begin{figure}[ht]
\begin{picture}(0,50)
\hspace{2cm}\includegraphics[scale=1]{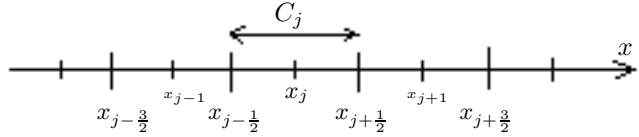}
\put(-150,43) {$C_j$}
\put(-146,14) {{\small $x_j$}}
\put(-127,6) {{\small $x_{j+\frac{1}{2}}$}}
\put(-175,6) {{\small $x_{j-\frac{1}{2}}$}}
\put(-217,6) {{\small $x_{j-\frac{3}{2}}$}}
\put(-80,6) {{\small $x_{j+\frac{3}{2}}$}}
\put(-192,14) {{\tiny $x_{j-1}$}}
\put(-100,14) {{\tiny $x_{j+1}$}}
\put(-20,30) {$x$}
\end{picture}
\caption{\textbf{Discr\'etisation spatiale}}
\end{figure}

\medskip
Pour la discr\'etisation temporelle, on se donne un pas de temps $\Delta t$ et une suite d'instants $t^n=n\Delta t,\; n\geq 0$.

\subsection{Un sch\'ema de type Godunov}

La topographie $a$ du probl\`eme (\ref{e2}) est approch\'ee par des valeurs discr\`etes $a_j,\;j\in J$ :

\begin{eqnarray}\label{e4} a_j\approx \dfrac{1}{\Delta x}\int_{C_j}a(x)dx.\end{eqnarray}

En posant :

\begin{eqnarray}\label{e14} a_\Delta(x)=\sum a_j\chi_j(x),\end{eqnarray}
avec $\chi_j(x)=1$ si $x\in C_j$ et $\chi_j(x)=0$ sinon, le syst\`eme (\ref{e2}) est approch\'e par : 

\begin{eqnarray}\label{e3}\begin{array}\{{lcl}.
\displaystyle\frac{\partial h}{\partial t}+\frac{\partial q}{\partial x}&=&0\\&&\\\displaystyle\frac{\partial q}{\partial t}+ \frac{\partial }{\partial x}\left(\dfrac{q^2}{h} + \frac{1}{2}gh^2 \right)+gh\dfrac{\partial a_\Delta}{\partial x}&=&0\\&&\\\dfrac{\partial a_\Delta}{\partial t}&=&0\end{array},(x,t)\in\R\times\R_+\end{eqnarray}

Posons $w=\begin{array}({c})h\\q\\a_\Delta\end{array}$. La solution $w$ du probl\`eme (\ref{e3}) est approch\'ee par des valeurs discr\`etes $w_j^n,j\in J, n\in \N$ :

\begin{eqnarray}\label{e15} w_j^n\approx \dfrac{1}{\Delta x}\int_{C_j}w(x,t^n)dx.\end{eqnarray}

\medskip
Puisque, sur chaque cellule $C_j$, $\dfrac{\partial a_\Delta}{\partial x}=0$, l'int\'egration de (\ref{e3}) sur $C_j\times [t^n,t^{n+1}]$ nous donne :

{\small $$\begin{array}{rcl}\displaystyle\int_{t^n}^{t^{n+1}}\int_{C_j}\left\{\frac{\partial w}{\partial t}+\frac{\partial F(w)}{\partial x}\right\}dxdt&=&0\\&&\\\displaystyle\int_{C_j}\left\{w(x,t^{n+1})-w(x,t^n)\right\}dx+\int_{t^n}^{t^{n+1}}\left\{F(w(x_{j+\frac{1}{2}}^-,t))-F(w(x_{j-\frac{1}{2}}^+,t))\right\}dt&=&0\end{array},$$} o\`u $F(w)=\begin{array}({c})q\\\frac{q^2}{h}+\frac{1}{2}gh^2\\0\end{array}$; $x^-$ et $x^+$ d\'esignant respectivement les limites \`a gauche et \`a droite de $x$.

\medskip
La premi\`ere int\'egrale de cette derni\`ere relation est approch\'ee en utilisant (\ref{e15}). Il restera alors l'approximation num\'erique de la deuxi\`eme int\'egrale.

\medskip
La m\'ethode des volumes finis repose sur le fait qu'\`a tout instant, la solution $w$ est constante par cellule. Ainsi, partant de la solution $w(x,t^n)$ \`a l'instant $t^n$, le calcul de $w(x_{j+\frac{1}{2}}^-,t)$ et $w(x_{j-\frac{1}{2}}^+,t)$ pour $t\in \left[t^n,t^{n+1}\right]$ est donn\'e par la r\'esolution du probl\`eme de Riemann suivant :

\begin{eqnarray}\label{e5}\begin{array}\{{lcl}.
\displaystyle\frac{\partial h}{\partial t}+\frac{\partial q}{\partial x}&=&0\\&&\\\displaystyle\frac{\partial q}{\partial t}+ \frac{\partial }{\partial x}\left(\dfrac{q^2}{h} + \frac{1}{2}gh^2 \right)+gh\dfrac{\partial a_\Delta}{\partial x}&=&0\\&&\\\dfrac{\partial a_\Delta}{\partial t}&=&0\\\\w(x,t^n)=\begin{array}\{{ll}.w_L,&\textrm{si }x<x_{j+\frac{1}{2}}\\w_R,&\textrm{si }x>x_{j+\frac{1}{2}}\end{array} \end{array},\end{eqnarray} o\`u $w_L=w_j^n$ et $w_R=w_{j+1}^n$.

\medskip
D\'esignons par $w_{j+1/2}^n(x/t;w_L,w_R)$ la solution autosimilaire de (\ref{e5}).

On d\'efinit deux flux $$F_{j+\frac{1}{2}}^{n,-}=F(w_{j+1/2}^n(0^-;w_L,w_R))\textrm{ et }F_{j+\frac{1}{2}}^{n,+}=F(w_{j+1/2}^n(0^+;w_L,w_R))$$ correspondant \`a chaque c\^ot\'e de l'interface de $x_{j+\frac{1}{2}}$. On obtient alors le sch\'ema num\'erique suivant :

\begin{eqnarray}\label{e6}w_j^{n+1}=w_j^n-\dfrac{\Delta t}{\Delta x}\left(F_{j+\frac{1}{2}}^{n,-}-F_{j+\frac{1}{2}}^{n,+}\right)\end{eqnarray}

\section{R\'esolution du probl\`eme de Riemann}

Il vient de ce qui pr\'ec\`ede que l'impl\'ementation num\'erique des \'equations de Saint Venant r\'eside en la r\'esolution, \`a chaque interface, du probl\`eme de Riemann (\ref{e5}); puisque c'est cette r\'esolution qui fournira les flux num\'eriques \`a utiliser.

\medskip
En posant $c=\sqrt{gh}$, \`a chaque interface $x=0$ et en dehors des zones s\`eches, (\ref{e5}) est \'equivalent \`a :

\begin{eqnarray}\label{e16}\begin{array}\{{lcl}.
\displaystyle\frac{\partial (2c)}{\partial t}+u\frac{\partial (2c)}{\partial x}+c\frac{\partial u}{\partial x}&=&0\\&&\\\dfrac{\partial u}{\partial t}+ c\dfrac{\partial (2c)}{\partial x}+u\dfrac{\partial u}{\partial x}+g\dfrac{\partial a_\Delta}{\partial x}&=&0\\&&\\\dfrac{\partial a_\Delta}{\partial t}&=&0\\\\w(x,0)=\begin{array}\{{ll}.w_L,&\textrm{si }x<0\\w_R,&\textrm{si }x>0\end{array} \end{array}\end{eqnarray}

\medskip
En posant $Y(w)=\begin{array}({c})2c\\u\\a_\Delta
\end{array}
$, on obtient finalement :

\begin{eqnarray}\label{e17}\begin{array}\{{l}.
\dfrac{\partial Y}{\partial t}+A(Y)\dfrac{\partial Y}{\partial x}=0\\\\Y(x,0)=\begin{array}\{{ll}.Y_L,&\textrm{si }x<0\\Y_R,&\textrm{si }x>0\end{array} \end{array}\end{eqnarray}
avec $A(Y)=\begin{array}({ccc})u&c&0\\c&u&g\\0&0&0
\end{array}$.

\subsection{Lin\'earisation du probl\`eme de Riemann}

A chaque interface, on r\'esout le syst\`eme lin\'eaire suivant :

\begin{eqnarray}\label{e7}\begin{array}\{{l}.
\displaystyle\frac{\partial Y}{\partial t}+A(\hat{Y})\frac{\partial Y}{\partial x}=0\\\\Y(x,0)=Y_0(x)=\begin{array}\{{l}.Y_L,\textrm{ si }x<0\\\\Y_R,\textrm{ si }x>0\end{array}\end{array}\end{eqnarray}
o\`u $\hat{Y}=\hat{Y}(Y_L,Y_R)$ est un \'etat moyen d\'ependant uniquement des \'etats $Y_L$ et $Y_R$; $\hat{Y}$ doit en plus v\'erifier la condition de consistance $\hat{Y}(Y^*,Y^*)=Y^*$. Un choix possible est donc : 
\begin{eqnarray}\label{e13}
\hat{Y}=\dfrac{Y_R+Y_L}{2}
\end{eqnarray}

\subsection{Solveur de Riemann}

$A(\hat{Y})$ admet trois valeurs propres qui sont :

\begin{eqnarray}\label{e12}
\hat{\lambda}_0=0, \hat{\lambda}_1=\hat{u}-\hat{c}, \hat{\lambda}_2=\hat{u}+\hat{c},
\end{eqnarray}

\medskip
Sous l'hypoth\`ese que $|\hat{u}|\neq \hat{c}$ et en dehors  des zones s\`eches, les trois valeurs propres sont deux \`a deux distinctes et $A(\hat{Y})$ est diagonisable.

\medskip
A chaque valeur propre $\hat{\lambda}_k$, nous associons un vecteur propre \`a droite $r_k\in\R^3$ : \begin{eqnarray}A(\hat{Y})r_k=\hat{\lambda}_kr_k\end{eqnarray} et un vecteur propre \`a gauche ${}^tl_k$ :\begin{eqnarray}{}^tl_kA(\hat{Y})=\hat{\lambda}_k{}^tl_k\end{eqnarray}

\begin{prop}\label{e10}${}$

$${}^tl_jr_k=0,\;\forall j\neq k.$$
\end{prop}

\subsubsection*{Preuve}

$\begin{array}{lll}{}^tl_jA=\hat{\lambda}_j{}^tl_j&\Longrightarrow&\left({}^tl_jA\right)r_k=\left(\hat{\lambda}_j{}^tl_j\right)r_k\\&&\\&\Longrightarrow&{}^tl_j\left(Ar_k\right)=\hat{\lambda}_j{}^tl_jr_k\\&&\\&\Longrightarrow&{}^tl_j(\hat{\lambda}_kr_k)=\hat{\lambda}_j{}^tl_jr_k\\&&\\&\Longrightarrow&\left(\hat{\lambda}_k-\hat{\lambda}_j\right){}^tl_jr_k=0\\&&\\&\Longrightarrow&{}^tl_jr_k=0,\;\forall j\neq k.\end{array}$
\begin{flushright}
{\tiny $\square$}
\end{flushright}

\begin{rmq}${}$

$\bullet$ ${}^tl_jr_k=l_j\cdot r_k.$

$\bullet$ Si $r$ est un vecteur propre de la valeur propre $\hat{\lambda}$, alors pour tout $\alpha\neq 0$, $\alpha r$ est \'egalement un vecteur propre de $\hat{\lambda}$. On peut donc choisir les bases ($r_0,r_1,r_2$) et ($l_0,l_1,l_2$) telles que $l_j\cdot r_j=1$, $\forall j=0,1,2.$
\end{rmq}

\medskip
On peut prendre :

\medskip
$\bullet$ $r_0=\begin{array}(c) \frac{g(\hat{\lambda}_2-\hat{\lambda}_1)}{2}\\\\-\frac{g(\hat{\lambda}_2+\hat{\lambda}_1)}{2}\\\\\hat{\lambda}_2\hat{\lambda}_1\end{array}$, $r_1=\begin{array}(c)-1\\1\\0\end{array}$,  $r_2=\begin{array}(c)1\\1\\0\end{array}$.

\medskip
et

\medskip
$\bullet$ $l_0=\begin{array}(c) 0\\\\0\\\\\frac{1}{\hat{\lambda}_2\hat{\lambda}_1}\end{array}$, $l_1=\begin{array}(c)-\frac{1}{2}\\\\\frac{1}{2}\\\\\frac{g}{2\hat{\lambda}_1}\end{array}$,  $l_2=\begin{array}(c)\frac{1}{2}\\\\\frac{1}{2}\\\\\frac{g}{2\hat{\lambda}_2}\end{array}$.

\medskip
(\ref{e7}) \'etant strictement hyperbolique, alors $(r_0,r_1,r_2)$ est une base de $\R^3$.

En d\'esignant par $Y^*$ la solution de (\ref{e7}), on peut alors \'ecrire : 

$Y^*(x,t)=\alpha_0(x,t)r_0+\alpha_1(x,t)r_1+\alpha_2(x,t)r_2$ et

$Y^*_0(x)=\alpha_0^0(x)r_0+\alpha_1^0(x)r_1+\alpha_2^0(x)r_2$

$$\begin{array}{rcl}\dfrac{\partial Y}{\partial t}+A(\hat{Y})\dfrac{\partial Y}{\partial x}=0&\Leftrightarrow&\displaystyle\sum_{i=0}^2\left\{\dfrac{\partial \alpha_i}{\partial t}r_i+\left(A(\hat{Y})r_i\right)\frac{\partial \alpha_i}{\partial x}\right\}=0\\&&\\&\Leftrightarrow&\displaystyle\sum_{i=0}^2\left\{\dfrac{\partial \alpha_i}{\partial t}r_i+\left(\hat{\lambda}_ir_i\right)\frac{\partial \alpha_i}{\partial x}\right\}=0\\&&\\&\Leftrightarrow&\displaystyle\sum_{i=0}^2\left\{\dfrac{\partial \alpha_i}{\partial t}+\hat{\lambda}_i\dfrac{\partial \alpha_i}{\partial x}\right\}r_i=0\\&&\\&\Leftrightarrow&\dfrac{\partial \alpha_i}{\partial t}+\hat{\lambda}_i\dfrac{\partial \alpha_i}{\partial x}=0, \forall i=0,1,2\end{array}$$

D'apr\`es la proposition\ref{e10}, on a :

\medskip
$\alpha_i^0(x)=\dfrac{1}{r_i\cdot l_i}Y^*_0(x)\cdot l_i=Y^*_0(x)\cdot l_i=\begin{array}\{{l}.\alpha_i^L=Y_L\cdot l_i,\textrm{ si }x<0\\\\\alpha_i^R=Y_R\cdot l_i,\textrm{ si }x>0\end{array}$

\medskip
Mais pour chaque $i$, l'\'equation d'advection \begin{eqnarray}\begin{array}\{{l}.\dfrac{\partial \alpha_i}{\partial t}+\hat{\lambda}_i\dfrac{\partial \alpha_i}{\partial x}=0\\\\\alpha_i^0(x)=\begin{array}\{{l}.\alpha_i^L,\textrm{ si }x<0\\\alpha_i^R,\textrm{ si }x>0\end{array}\end{array}\end{eqnarray} admet pour unique solution : $\alpha_i(x,t)=\begin{array}\{{l}.\alpha_i^L,\textrm{ si }\frac{x}{t}<\hat{\lambda}_i\\\\\alpha_i^R,\textrm{ si }\frac{x}{t}>\hat{\lambda}_i\end{array}$

Ainsi la solution du probl\`eme de Riemann (\ref{e7}) est :

$$\begin{array}{lll}
Y^*(x,t)&=&\displaystyle Y_L+\sum_{x/t>\hat{\lambda}_k}(\alpha_k^R-\alpha_k^L)r_k\\&&\\&=&\displaystyle Y_R-\sum_{x/t<\hat{\lambda}_k}(\alpha_k^R-\alpha_k^L)r_k
  \end{array}$$

\subsection{Solutions du probl\`eme de Riemann}\label{e18}

Nous pr\'esentons dans ce paragraphe les trois cas possibles de solution du probl\`eme de Riemann (\ref{e7}) dans le plan $(x,t)$.
On d\'esigne $Y^*(0^-,t)$ et $Y^*(O^+,t)$ par $Y^*_l$ et $Y^*_r$ respectivement.

\begin{figure}[ht]
\begin{picture}(0,90)
\hspace{2cm}\includegraphics[scale=.8]{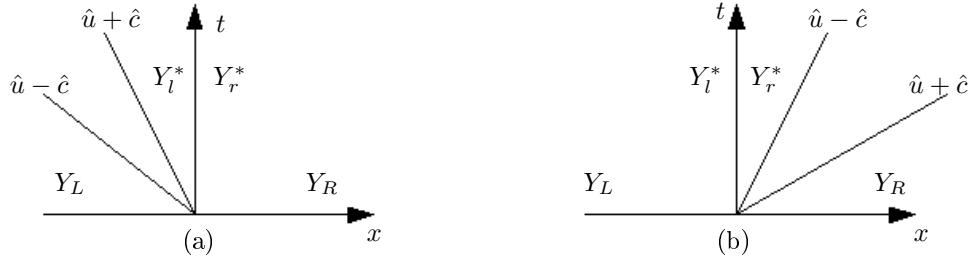}
\put(-222,1) {$x$}
\put(-279,80) {$t$}
\put(-340,20) {$Y_L$}
\put(-245,20) {$Y_R$}
\put(-303,60) {$Y^*_l$}
\put(-280,60) {$Y^*_r$}
\put(-357,57) {$\hat{u}-\hat{c}$}
\put(-330,82) {$\hat{u}+\hat{c}$}
\put(-292,-2) {(a)}

\put(-15,1) {$x$}
\put(-90,85) {$t$}
\put(-140,20) {$Y_L$}
\put(-30,20) {$Y_R$}
\put(-100,60) {$Y^*_l$}
\put(-77,60) {$Y^*_r$}
\put(-17,57) {$\hat{u}+\hat{c}$}
\put(-55,82) {$\hat{u}-\hat{c}$}
\put(-90,-2) {(b)}
\end{picture}

\caption{\textbf{Ecoulement torrentiel}}
\end{figure}

\vspace{-.5cm}
$$\textrm{(a)}:\begin{array}\{{l}.Y^*_l=Y^*_r-(\alpha_0^R-\alpha_0^L)r_0\\\\Y^*_r=Y_R
\end{array}; \textrm{(b)}:\begin{array}\{{l}.Y^*_l=Y_L\\\\Y^*_r=Y^*_l+(\alpha_0^R-\alpha_0^L)r_0
\end{array}$$

\begin{figure}[ht]
\begin{picture}(0,100)
\hspace{2cm}\includegraphics[scale=.8]{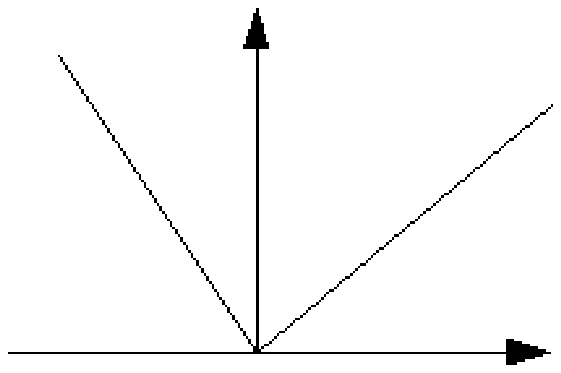}
\put(-15,4) {$x$}
\put(-85,89) {$t$}
\put(-130,23) {$Y_L$}
\put(-25,23) {$Y_R$}
\put(-92,67) {$Y^*_l$}
\put(-68,67) {$Y^*_r$}
\put(-17,77) {$\hat{u}+\hat{c}$}
\put(-136,85) {$\hat{u}-\hat{c}$}
\put(40,50) {$\begin{array}\{{l}.Y^*_l=Y_L+(\alpha_1^R-\alpha_1^L)r_1\\\\Y^*_r=Y_R-(\alpha_2^R-\alpha_2^L)r_2
\end{array}$}
\end{picture}

\caption{\textbf{Ecoulement fluvial}}\label{e11}
\end{figure}

\newpage
\subsection{Correction entropique}

Le solveur de Roe peut conduire \`a une c\'el\'erit\'e n\'egative faisant ainsi exploser l'impl\'ementation num\'erique.

Par exemple, prenons un \'ecoulement sur fond plat ($a(x)\equiv 0$) avec $u_L=-\dfrac{3}{2}\hat{c}$, $u_R=\dfrac{5}{2}\hat{c}$.

On a alors ((\ref{e13}) et (\ref{e12})): $\hat{\lambda}_1=-\dfrac{1}{2}\hat{c}<0$; $\hat{\lambda}_2=\dfrac{3}{2}\hat{c}>0$.

L'\'ecoulement local est donc fluvial. D\'eterminons la composante $2c^*_l$ de $Y^*_l$ (voir FIG.\ref{e11}) : 

$$\begin{array}{lll}
2c^*_l&=&2c_L+(\alpha_1^R-\alpha_1^L)r_1^1\quad \textrm{ (avec $r_1^1$, premi\`ere composante de $r_1$)}\\&&\\&=&2c_L-\left(-\dfrac{1}{2}(2c_R-2c_L)+\dfrac{1}{2}(u_R-u_L)\right)\\&&\\&=&2c_L-\left(-\dfrac{1}{2}(2c_R-2c_L)+\dfrac{1}{2}\times 4\times\dfrac{2c_R+2c_L}{2}\right)\\&&\\&=&-(c_R+c_L)<0
\end{array}
$$

\subsubsection{C\'el\'erit\'e de l'\'etat interm\'ediaire}

En reprenant les trois configurations de la section-\ref{e18}, on obtient :

\textbf{\uline{Cas (a)}} : 

\begin{eqnarray}\label{e19}\begin{array}\{{l}. 2c^*_l=2c_R-\dfrac{g(a_R-a_L)(\hat{\lambda}_2-\hat{\lambda}_1)}{2\hat{\lambda}_2\hat{\lambda}_1}\\\\2c^*_r=2c_R\end{array}
\end{eqnarray}

\textbf{\uline{Cas (b)}} : 

\begin{eqnarray}\label{e20}\begin{array}\{{l}.2c^*_l=2c_L\\\\2c^*_r=2c_L+\dfrac{g(a_R-a_L)(\hat{\lambda}_2-\hat{\lambda}_1)}{2\hat{\lambda}_2\hat{\lambda}_1}\end{array}
\end{eqnarray}

\textbf{\uline{Cas (c)}} : 

\begin{eqnarray}\label{e21}\begin{array}\{{l}.2c^*_l=\dfrac{1}{2}\left(2(c_R+c_L)-(u_R-u_L)-\dfrac{g(a_R-a_L)}{\hat{\lambda}_1}\right)\\\\2c^*_r=\dfrac{1}{2}\left(2(c_R+c_L)-(u_R-u_L)-\dfrac{g(a_R-a_L)}{\hat{\lambda}_2}\right)\end{array}
\end{eqnarray}

\subsubsection{Ecoulement sur fond plat}

\begin{prop}\cite{gallouet3}
Si \begin{eqnarray}\label{e22}
u_R-u_L<2(c_R+c_L)
\end{eqnarray}
alors l'\'etat interm\'ediaire $Y^*$ admet une c\'el\'erit\'e positive.
\end{prop}

\subsubsection*{Preuve}

La preuve est imm\'ediate; en effet, ici $a_R-a_L=0$. Par suite, la c\'el\'erit\'e de l'\'etat interm\'ediaire $Y^*$ est positive d'apr\`es \ref{e19},\ref{e20} et \ref{e21}.
\begin{flushright}
{\tiny $\square$}
\end{flushright}

Si la condition (\ref{e22}) est viol\'ee, on pose \cite{gallouet3} : $c^*_r=c^*_l=0$.
\subsubsection{Ecoulement avec topographie}

\medskip
On choisit alors d'agir sur les vitesses d'ondes $\hat{\lambda}$ de l'\'etat interm\'ediaire $Y^*$; on notera $\lambda$ les nouvelles vitesses de $Y^*$ qui seront utilis\'ees dans l'impl\'ementation num\'erique.

\begin{prop}
On suppose $u_R-u_L<2(c_R+c_L)$;  si $\hat{\lambda}_1<0<\hat{\lambda}_2$, le choix $$\begin{array}\{{l}.\lambda_1=\hat{\lambda}_1\\\lambda_2=max\left(\hat{\lambda}_2,\dfrac{g(a_R-a_L)}{2(c_R+c_L)-(u_R-u_L)}\right)
\end{array},\; \textrm{si } a_R>a_L$$ et $$\begin{array}\{{l}.\lambda_1=min\left(\hat{\lambda}_1,\dfrac{g(a_R-a_L)}{2(c_R+c_L)-(u_R-u_L)}\right)\\\lambda_2=\hat{\lambda}_2
\end{array},\; \textrm{si } a_R<a_L
$$ permet l'obtention d'un \'etat interm\'ediaire $Y^*$ \`a c\'el\'erit\'e positive.
\end{prop}

\subsubsection*{Preuve}

La preuve est imm\'ediate d'apr\`es \ref{e21}.
\begin{flushright}
{\tiny $\square$}
\end{flushright}

Si la condition (\ref{e22}) est viol\'ee, on pose : 

$\begin{array}\{{l}. c^*_l=-\dfrac{g(a_R-a_L)}{4\hat{\lambda}_1}\\c^*_r=0
\end{array}$, si $a_R>a_L$ et $\begin{array}\{{l}. c^*_l=0\\c^*_r=-\dfrac{g(a_R-a_L)}{4\hat{\lambda}_2}
\end{array}$, si $a_R<a_L$

\begin{prop}
Si $\hat{\lambda}_2<0$, le choix $$\begin{array}\{{l}.\lambda_1=\hat{\lambda}_1\\\lambda_2=min\left(\hat{\lambda}_2,-\dfrac{g(a_R-a_L)\hat{\lambda}_1}{4c_R\hat{\lambda}_1-g(a_R-a_L)}\right)
\end{array},\; \textrm{si } a_R>a_L \textrm{ et }\begin{array}\{{l}.\lambda_1=\hat{\lambda}_1\\\lambda_2=\hat{\lambda}_2
\end{array},\; \textrm{si } a_R<a_L
$$ permet l'obtention d'un \'etat interm\'ediaire $Y^*$ \`a c\'el\'erit\'e positive.
\end{prop}

\subsubsection*{Preuve}

La preuve est imm\'ediate d'apr\`es \ref{e19}.
\begin{flushright}
{\tiny $\square$}
\end{flushright}

\begin{prop}
Si $\hat{\lambda}_1>0$, le choix $$\begin{array}\{{l}.\lambda_1=\hat{\lambda}_1\\\lambda_2=\hat{\lambda}_2
\end{array},\; \textrm{si } a_R>a_L \textrm{ et }\begin{array}\{{l}.\lambda_1=max\left(\hat{\lambda}_1,-\dfrac{g(a_R-a_L)\hat{\lambda}_2}{4c_L\hat{\lambda}_2-g(a_R-a_L)}\right)\\\lambda_2=\hat{\lambda}_2
\end{array},\; \textrm{si } a_R<a_L$$ permet l'obtention d'un \'etat interm\'ediaire $Y^*$ \`a c\'el\'erit\'e positive.
\end{prop}

\subsubsection*{Preuve}

La preuve est imm\'ediate d'apr\`es \ref{e20}.
\begin{flushright}
{\tiny $\square$}
\end{flushright}

\subsection{Traitement des bancs-couvrants-d\'ecouvrants}
On se situe ici dans le cas o\`u la hauteur de l'\'etat $w_R$ ou celle de $w_L$ est nulle. Pour la d\'etermination de l'\'etat interm\'ediaire $Y^*$, nous utilisons les r\'esultats de \cite{david} bas\'es sur une nouvelle estimation des vitesses d'onde $\lambda_1$ et $\lambda_2$ : Une onde de d\'etente est g\'en\'er\'ee du c\^ot\'e o\`u la hauteur d'eau est nulle \cite{david, leveque}.

\begin{figure}[ht]
\begin{picture}(0,160)
\hspace{2.2cm}\includegraphics[scale=.8]{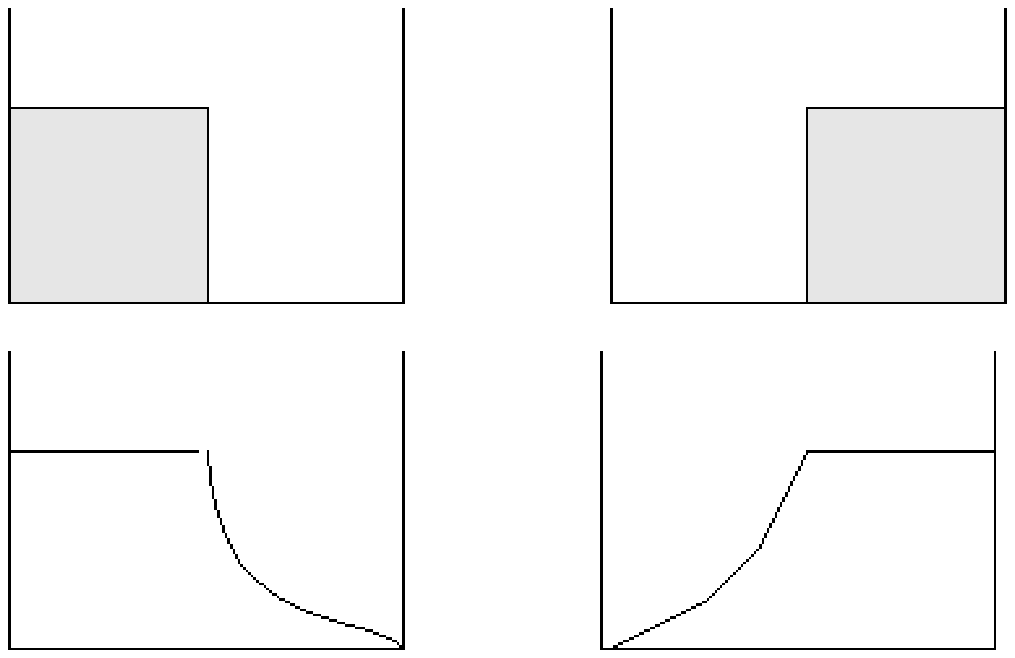}
\put(-52,120) {$h_R>0$}
\put(-96,91) {$h_L=0$}
\put(-235,120) {$h_L>0$}
\put(-185,91) {$h_R=0$}
\put(-235,-17) {$\begin{array}\{{l}.\lambda_1=u_L-c_L\\\lambda_2=u_L+2c_L
\end{array}
$}
\put(-95,-17) {$\begin{array}\{{l}.\lambda_1=u_L-2c_L\\\lambda_2=u_L+c_L
\end{array}
$}
\end{picture}

\vspace*{.7cm}
\caption{bancs-couvrants-d\'ecouvrants}
\end{figure}

\newpage
\section{R\'esultats num\'eriques}

\subsection{Ecoulements noy\'es}

La longueur du domaine de calcul est $L=25m$, le nombre de points du maillage $N=1000$; la dur\'ee totale d'observation est $T=1.2s$; les conditions initiales sont :

\begin{center}
$\begin{array}{|c|c|c|}\hline &0\leq x\leq 12.5&12.5\leq x\leq 25\\\hline h(x,t=0)&3&4\\\hline u(x,t=0)&0&0\\\hline
\end{array}$
\end{center}

\newpage
\subsubsection{Test1 : Ecoulement sur fond plat}
$a(x)=0,\forall x\in [0;25]$; $CFL=0.8$.

\vspace{.55cm}
\begin{figure}[ht]
\begin{picture}(0,170)
\hspace{1.5cm}\includegraphics[scale=.5]{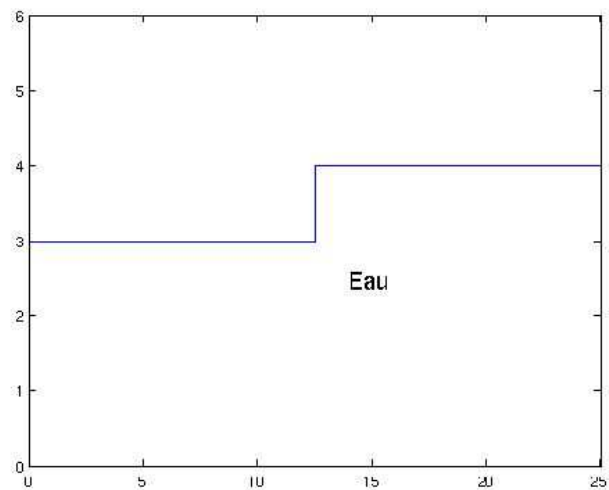}
\end{picture}
\caption{Test1 (t=0s)}
\end{figure}

\vspace{1cm}
\begin{figure}[ht]
\begin{picture}(0,170)
\hspace{1.5cm}\includegraphics[scale=.5]{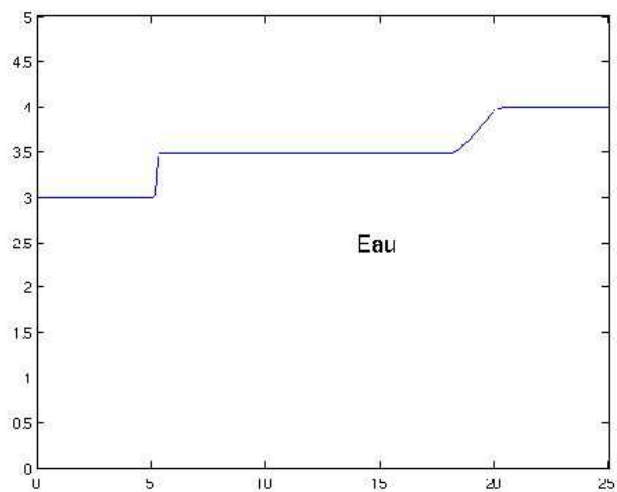}
\end{picture}
\caption{Test1 (t=1.2s)}
\end{figure}

\newpage
\subsubsection{Test2 : Ecoulement avec topographie}
$a(x)=2$ sur $[0;12.5]$ et $a(x)=0$ sur $[12.5;25]$; $CFL=0.8$.

\vspace{.55cm}
\begin{figure}[ht]
\begin{picture}(0,170)
\hspace{1.5cm}\includegraphics[scale=.5]{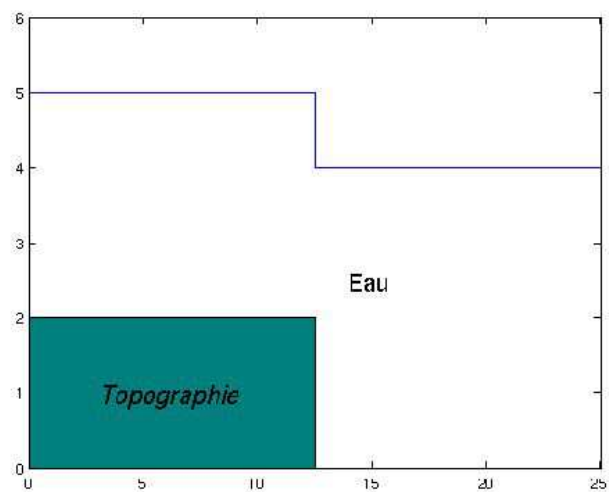}
\end{picture}
\caption{Test2 (t=0s)}
\end{figure}

\vspace{1cm}
\begin{figure}[ht]
\begin{picture}(0,170)
\hspace{1.5cm}\includegraphics[scale=.5]{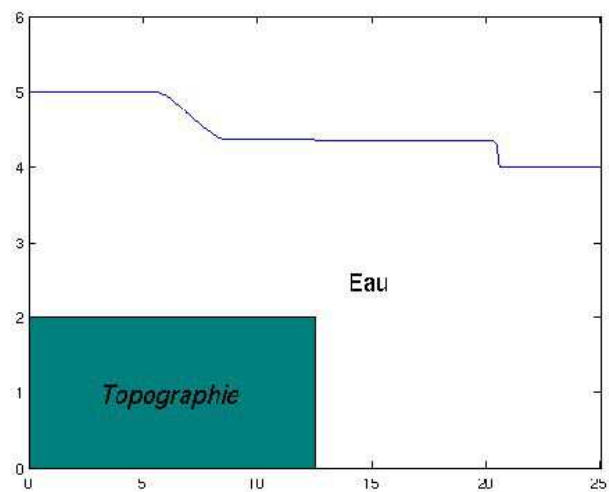}
\end{picture}
\caption{Test2 (t=1.2s)}
\end{figure}

\newpage
\subsubsection{Test3 : Ecoulement avec topographie}
$a(x)=4$ sur $[0;12.5]$ et $a(x)=0$ sur $[12.5;25]$; $CFL=0.6$.

\vspace{.55cm}
\begin{figure}[ht]
\begin{picture}(0,170)
\hspace{1.5cm}\includegraphics[scale=.5]{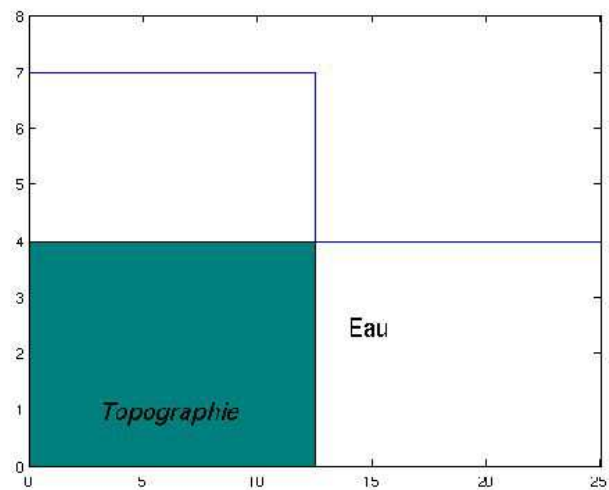}
\end{picture}
\caption{Test3 (t=0s)}
\end{figure}

\vspace{1cm}
\begin{figure}[ht]
\begin{picture}(0,170)
\hspace{1.5cm}\includegraphics[scale=.5]{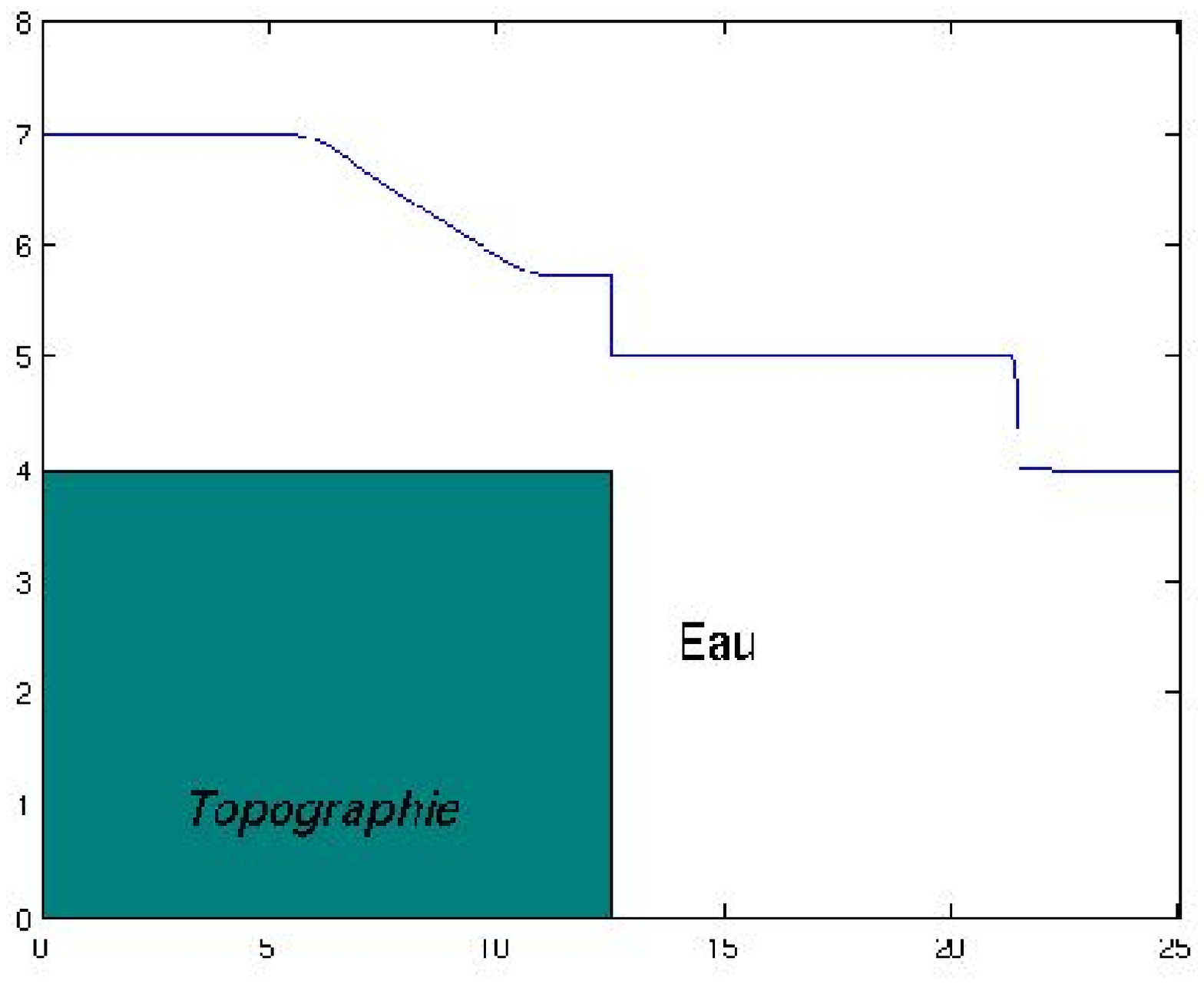}
\end{picture}
\caption{Test3 (t=1.2s)}
\end{figure}

\newpage
\subsection{Bancs-couvrants-d\'ecouvrants}

La longueur du domaine de calcul est $L=25m$, le nombre de points de maillage $N=500$; la dur\'ee totale d'observation est $T$; les conditions initiales sont :

\subsubsection*{Test 4}

\begin{center}
$\begin{array}{|c|c|c|}\hline &0\leq x\leq 12.5&12.5\leq x\leq 25\\\hline h(x,t=0)&1.5&0\\\hline u(x,t=0)&0&0\\\hline
\end{array}$
\end{center}

\subsubsection*{Test 5}

\begin{center}
$\begin{array}{|c|c|c|c|c|}\hline &0\leq x\leq 7.5&7.5\leq x\leq 10.5&10.5\leq x\leq 14.5&14.5\leq x\leq 25\\\hline h(x,t=0)&1&0.1-a(x)&0&0.1-a(x)\\\hline u(x,t=0)&0&0&0&0\\\hline
\end{array}$
\end{center}
\newpage
\subsubsection{Test4 : Bancs-couvrants-d\'ecouvrants}
$a(x)=0$ sur $[0;12.5]$ et $a(x)=1$ sur $[12.5;25]$;  $T=2s$; $CFL=0.4$.

\vspace{.55cm}
\begin{figure}[ht]
\begin{picture}(0,170)
\hspace{1.5cm}\includegraphics[scale=.5]{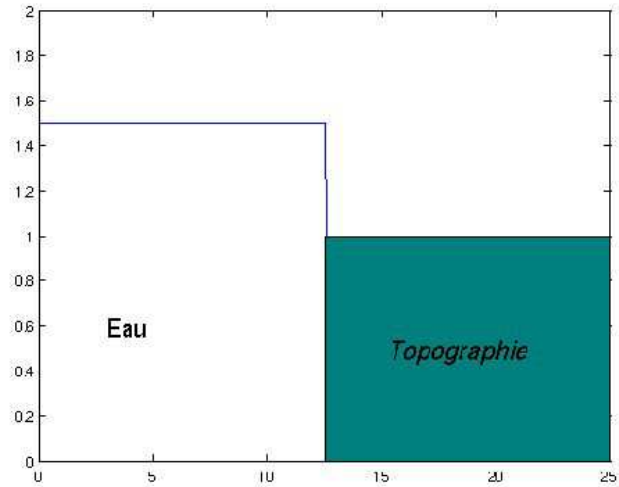}
\end{picture}
\caption{Test4 (t=0s)}
\end{figure}

\vspace{1cm}
\begin{figure}[ht]
\begin{picture}(0,170)
\hspace{1.5cm}\includegraphics[scale=.5]{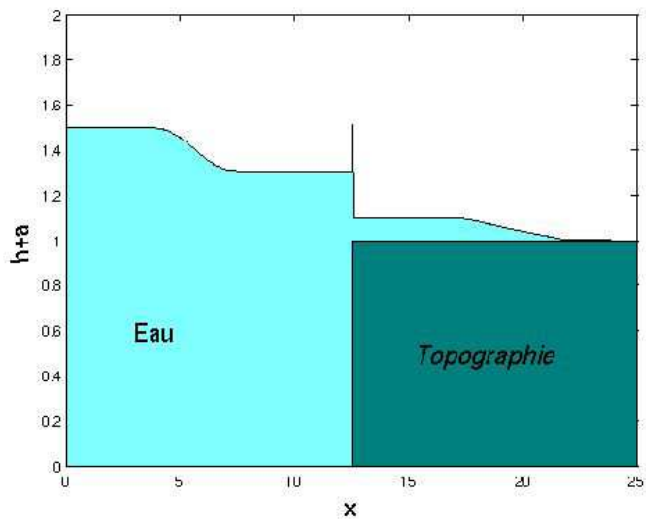}
\end{picture}
\caption{Test4 (t=2s)}
\end{figure}

\newpage
\subsubsection{Test5 : Bancs-couvrants-d\'ecouvrants}
$a(x)$ (voir figure);  $T=3s$; $CFL=0.4$.

\vspace{.55cm}
\begin{figure}[ht]
\begin{picture}(0,170)
\hspace{1.5cm}\includegraphics[scale=.5]{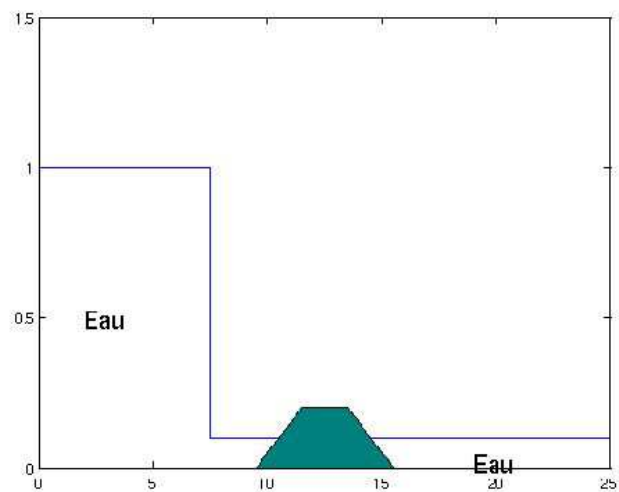}
\end{picture}
\caption{Test5 (t=0s)}
\end{figure}

\vspace{1cm}
\begin{figure}[ht]
\begin{picture}(0,170)
\hspace{1.5cm}\includegraphics[scale=.5]{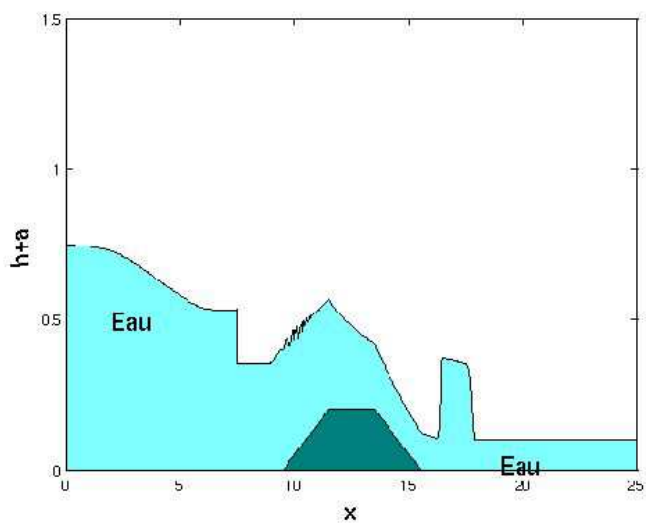}
\end{picture}
\caption{Test5 (t=3s)}
\end{figure}

\section{Conclusion}
Exploitant les r\'esultats de \cite{gallouet, gallouet2, gallouet3}, nous avons pr\'esent\'e une formulation c\'el\'erit\'e-vitesse des \'equations de Saint-Venant. La m\'ethode VFRoe conduit \`a la r\'esolution, par interface du maillage spatial, d'un probl\`eme de Riemann. Notre contribution a \'et\'e, en agissant sur les vitesses d'ondes, la construction d'un solveur de Riemann approch\'e garantissant la pr\'eservation de la positivit\'e de la c\'el\'erit\'e de l'\'etat interm\'ediaire dans la r\'esolution du probl\`eme de Riemann. L'exploitation des r\'esultats de \cite{david, leveque} nous a ensuite permis d'adapter le solveur afin de prendre en compte la simulation num\'erique des bancs-couvrants-d\'ecouvrants. 

Le solveur de Riemann approch\'e pr\'esent\'e dans ce rapport permet d'effectuer des simulations num\'eriques des \'equations de Saint-Venant, m\^eme en pr\'esence de hauteurs d'eau nulles (recouvrement, d\'ecouvrement). La prochaine \'etape de notre travail consistera \`a l'extension \`a des tests bidimensionnels.
\nocite{*}

\bibliographystyle{unsrt}
\newpage
\bibliography{mabiblio}

\end{document}